\documentclass[a4paper,notitlepage, 12pt,reqno]{article}

  \usepackage{mathtext}
 \usepackage[cp1251]{inputenc}
  \usepackage[T2A]{fontenc}
 \usepackage[russian,english]{babel}
 \usepackage[all,arc,poly,2cell,curve,arrow,tips]{xy}

  \usepackage{enumerate, eucal, amsthm,amsmath, amssymb}

  \allowdisplaybreaks[4]

 \theoremstyle{definition}
 \newtheorem{defn}{Определение}

 \theoremstyle{plain}
 \newtheorem{thm}{Tеорема}
 \newtheorem*{thm*}{Tеорема}
 \newtheorem{prop}{Предложение}
  \newtheorem*{prop*}{Предложение}
 
  \newtheorem*{cor*}{Следствие}
 \newtheorem{lem}{Лемма}
  \newtheorem*{lem*}{Лемма}

 \theoremstyle{remark}

 \newtheorem*{remark*}{Замечание}

 \renewcommand{\abstractname}{}

  \newcounter{ab}
\title{ Knizhnik-Zamolodchikov  type equations  for the root system $B$ and Capelli central elements}

 \author{D.V. Artamonov, V.A. Goloubeva}
  \date{}

\begin{document}

 \maketitle

\renewcommand{\abstractname}{}

\begin{abstract} 
 The
construction  of the well-known  Knizhnik-Zamolodchikov
  equations uses the central element of the second order in the universal enveloping  algebra for some
   Lie algebra. But in the universal enveloping algebra there are central elements of higher orders. It seems desirable to use these elements for the construction of Knizhnik-Zamolodchikov type 
  equations. In the present paper
    we give a construction of such   Knizhnik-Zamolodchikov type
  equations for the root system $B$ associated with Capelli central elements in the universal enveloping  algebra for the orthogonal algebra.
\end{abstract}

\maketitle{}

 \section{Introduction}

 The Knizhnik-Zamolodchikov equations are a
system of differential equations which is satisfied by correlation
functions in the WZW theory \cite{KZ}. Later it turned
out that these equation are related with many other areas of
mathematics (quantum algebra, isomonodromic
deformation).

The Kniznik-Zamolodchikov equations are also interesting  as a
nontrivial example of an integrable Pfaffian system of Fuchsian
type. Mention that the monodromy  representation of this system is
known explicitly. Thus we get a solution for the Riemann-Hilbert
problem in a very particular case.

Let us write the Kniznik-Zamolodchikov equations. They have the
form.

$$dy=\lambda(\sum_{i\neq
j=1}^n\tau_{ij}\frac{d(z_i-z_j)}{z_i-z_j})y,$$

where $\lambda$ is some  complex parameter, $y(z_1,...,z_n)$ is a
 vector function that takes values in a  tensor power
  $V^{\otimes n}$ of a representation space $V$ of a Lie algebra $\mathfrak{g}$, and
$\tau_{ij}$ is defined by the formula

 \begin{equation}\label{cladi}\tau_{ij}=\sum_s 1\otimes...\otimes \rho(I_{s})\otimes ...\otimes \rho(\omega(I_s))\otimes...1.\end{equation}

Here  $I_s$ is base of a finite-dimensional Lie algebra,
$\{\omega(I_s)\}$ is a base dual to $\{I_s\}$ with respect to the
Killing form. The elements $I_s$ occur on tensor factors $i$ and
$j$, and $\rho: \mathfrak{g}\rightarrow End(V)$ is a representation
of the Lie algebra $\mathfrak{g}$. Thus $\tau_{ij}$ is a matrix of
the operator on $V^{\otimes n}$.

The Knizhnik-Zamolodchikov equations have singularities on the
planes $z_i=z_j$. It is natural to look for similar systems that
have other singular locus.

Several authors constructed  different  Knizhnik-Zamolodchikov
 type equations, that have  singularities on the reflection hyperplanes
  of different  root system. Thus in  \cite{Ma}  a system  the 
     Knizhnik-Zamolodchikov  type equations  is constructed that have as  the singular set the reflection hyperplanes
     corresponding to an arbitrary root system.  Let us present such a  system.

Matsuo considers a root system $\Delta$, the set of positive roots
is denoted as $\Delta^+$ and the corresponding Weyl group  is
denoted as $W$. The function $y$ takes values in the group algebra
$\mathbb{C}[W]$ and the independent variable belong to the dual to
the vector space spanned by roots.

 Denote a reflection  corresponding to the root $\alpha$ as
 $\sigma_{\alpha}$.

 The Matsuo's  Knizhnik-Zamolodchikov type  equations contain arbitrary parameters $\lambda_{|\alpha|}$
  depending on the length of the roots  $\alpha$ and one additional parameter
  $\lambda$. The system has the following form.

$$\frac{\partial y}{\partial \xi}=(\sum_{\alpha\in \Delta^+}\lambda_{|\alpha|}(\alpha,z)\frac{1}{(\alpha,u)}(\sigma_{\alpha}-1))y.$$

Note that the solutions of this system take values in the  space
$\mathbb{C}[W]$
 and no Lie algebra is involved in it's construction explicitly.  Nevertheless if one takes in
the Matsuo's construction the root system $A$ and some special
representation of the Weyl group one can obtain from Matsuo's system
an ordinary Knizhnik-Zamolodchikov system corresponding to the
algebra $\mathfrak{sl}_n$ and its standard representation. Later
systems of Matsuo's type were intensively studied by numerous
authors (I. Cherednik \cite{ch} and others)

There exist other constructions of Knizhnik-Zamolodchikov type
equations with are closer to the original equations.  One of them is the
Leibman's construction of Knizhnik-Zamolodchikov type equations associated with the root
system $B$ \cite{L} and Enriquez's cyclotomic Knizhnik-Zamolodchikov type equations \cite{En}.

In \cite{AG} we have constructed  Knizhnik-Zamolodchikov type
equations  for the root system  $A$, but our construction was based on special centra elements in the universal enveloping algebra of the orthogoanl algebra, namely the Capelli elements. In the present paper we generalize our construction to the root system $B$ and show that  the construction  from \cite{AG} can be in fact  interpreted in some sence as a very special case  of the general Leibman's construction.

\section{Knizhnik-Zamolodchikov equation associated with the root system $B$}
\label{K}

Let us explain the Leiman's construction of   Knizhnik-Zamolodchikov type
equations  for the root system  $B$.

This is a system of type

$$dy = \lambda(\sum_{i<j}\frac{d(z_i-z_j)}{z_i-z_j}\tau_{ij}+\sum_{i<j}\frac{d(z_i-z_j)}{z_i-z_j}\mu_{ij}+\sum_i\frac{\nu_i}{z_i})y.$$

Here $\tau_{ij}$ is defined as in \ref{cladi}

 \begin{equation}\label{cladi}\tau_{ij}=\sum_s 1\otimes...\otimes \rho(I_{s})\otimes ...\otimes \rho(\omega(I_s))\otimes...1,\end{equation}

where  as before $I_s$ is a base of a Lie algebra,  $\{\omega(I_s)\}$ is a
base dual to $\{I_s\}$  with respect to the Killing form.

the coefficients $\mu_{ij}$ and $\nu_i$ are defined as follows.
 Denote as $\sigma$ an involution in the considered Lie algebra $\mathfrak{g}$. Then one has

 \begin{equation}\label{clam}\mu_{ij}=\sum_s 1\otimes...\otimes \rho(I_{s})\otimes ...\otimes \rho(\sigma(\omega(I_s)))\otimes...1,\end{equation}

 the elements $I_s$ occur on places $i,j$,

\begin{equation}\label{clan}\nu_{i}=\frac{1}{2}\sum_s 1\otimes.. ...\otimes \rho(I_s\sigma(\omega(I_s))+I_s^2)\otimes...1,\end{equation}

 the elements $I_s$ occur on  the place $i$.

 Although the Leibman's proof     in \cite{L} is done for the case of a simple Lie algebra and is  based on calculations in the
 root base in \cite{En} these is a proof that is valid for an
 arbitrary finite-dimensional Lie algebra with a fixed central element in $U(\mathfrak{g})$ of the second order.

\section{The Lie algebra $\mathcal{T}$.}

In this section we introduce a Lie algebra $\mathcal{T}$ which  playes the crucial role in our construction.

Consider the space of skew-symmetric tensors with $2n$ indices. Let
each index of the skew symmetric tensor take values in the set
$-n,...,n$.

There exist a structure of an associative algebra algebra on this
space

$$(e_{a_1}\wedge...\wedge e_{a_n}\wedge e_{b_1}\wedge...\wedge e_{b_n})\cdot (e_{-b_1}\wedge...\wedge e_{-b_n}\wedge e_{c_1}\wedge...\wedge e_{c_n}):=
e_{a_1}\wedge...\wedge e_{a_n}\wedge e_{c_1}\wedge...\wedge
e_{c_n}.$$

As a corollary we have a structure of a Lie algebra. Denote this algebra as $\mathcal{T}$.

This algebra has a representation on the space of skew-symmetric
tensors with $n$ indices defined by the formula

\begin{equation}\label{repr}e_{a_1}\wedge...\wedge e_{a_n}\wedge
e_{b_1}\wedge...\wedge e_{b_n} (e_{-b_1}\wedge...\wedge e_{-b_n}):=
 e_{a_1}\wedge...\wedge
e_{a_n}.\end{equation}

 The algebra has an involution $\omega$ defined as
follows

$$\omega(e_{a_1}\wedge...\wedge e_{a_n}\wedge e_{b_1}\wedge...\wedge e_{b_n})= e_{-a_1}\wedge...\wedge e_{-a_n}\wedge e_{-b_1}\wedge...\wedge e_{-b_n}.$$

In $U(\mathcal{T})$ there is a central element of the second order,
namely the element

$$C=\sum_{a_1,...,a_n,b_1,...,b_n}(e_{a_1}\wedge..\wedge e_{a_n}\wedge e_{b_1}\wedge...\wedge e_{b_n})\bullet (e_{-a_1}\wedge..\wedge e_{-a_n}\wedge e_{-b_1}\wedge...\wedge e_{-b_n}),$$

where $\bullet$ denotes the multiplication in $U(\mathcal{T})$.  The proof of this fact is essentialy  contained in \cite{Molev}.

Using general constructions described in  Section \ref{K} one can
construct  Knizhnik-Zamolodchikov type equations associated with
the root system $B$ based with coefficients in the Lie algebra
$\mathcal{T}$.

In the next section we give an interpretation of these equations as
Knizhnik-Zamolodchikov type equations whose construction is based on
some certain higher order central elements in
$U(\mathfrak{o}_{2n+1})$.

\section{Capelli elements and noncommutative pfaffians}
Let us define some certain central elements in the universal
enveloping algebra of the orthogonal algebra.

\subsection{The split realization of the orthogonal algebra}
We use the split realization of the orthogonal algebra. This means
that we define the orthogonal algebra as the algebra that preserves
the quadratic form

 \begin{align*}
G=&\begin{pmatrix} 0 & 0  &... & 0 &1\\
0& 0 & ...& 1 & 0\\
...\\
0 & 1 & ... & 0 & 0\\
1 & 0 &...& 0 & 0\\
\end{pmatrix}
  \end{align*}

The row and columns are indexed by $i,j=-n,-n+1,...,n-1,n$.  The
zero is skipped in the case $N=2n$ and is included in the case
$N=2n+1$. The algebra $\mathfrak{o}_N$ is generated by  matrices
$$F_{ij}=E_{ij}-E_{-j-i}.$$

The commutation relations between these generators are the following

\begin{equation}[F_{ij},F_{kl}]=\delta_{kj}F_{il}-\delta_{il}F_{kj}+\delta_{i-k}F_{-jl}+
\delta_{-lj}F_{k-i}.\end{equation}

\subsection{Noncommutative pfaffians and Capelli elements}

Now let us describe some special higher order central elements in the universal enveloping for the orthogoanl algebra.

Let  $\Phi=(\Phi_{ij})$ be a $k\times k$ matrix, where  $k$  is
even, whose elements belong to some noncommutative ring. The
noncommutative pfaffian is defined as follows:

$$Pf\Phi=\frac{1}{2^{\frac{k}{2}}(\frac{k}{2})!}\sum_{\sigma\in S_{k}}(-1)^{\sigma}\Phi_{\sigma(1)\sigma(2)}...\Phi_{\sigma(k-1)\sigma(k)}.$$

For a subset $I\subset \{-n,...,n\}$ define a submatrix
$F_I=(F_{ij})_{i,j\in I}$. For this subset put

$$Pf F_I:=Pf(F_{-i,j})_{-i,j\in I}.$$

\begin{defn}\label{P} Put

$C_k=\sum_{I\subset \{1,...,N\}, \,|I|=k}PfF_IPfF_I,$

The elements $C_k$ are the Capelli elements.

\end{defn}

\begin{thm}\cite{Molev}For odd $N$ the elements $C_k$ are algebraically independent and generate the center,
 for  even $N$ the same is true if one takes instead the highest Capelli element $C_N=(PfF)^2$ the central element  $PfF$.
\end{thm}

Below we need two formulas. There proofs can be found in \cite{Msb}.

\begin{lem}
\label{minors1}

$PfF_I=\frac{(\frac{p}{2})!(\frac{q}{2})!}{(\frac{k}{2})!}\sum_{I=I'\sqcup
I'',|I'|=p,|I''|=q}(-1)^{(I'I'')}PfF_{I'}PfF_{I''}$.

Here $(-1)^{(I'I'')}$ is a sign of a permutation of the set
$I=\{i_1,...,i_k\}$ that places first the subset $I'\subset I$ and
then the subset $I''\subset I$.

The numbers $p$, $q$ are even fixed numbers,  they satisfy
$p+q=k=|I|$.
\end{lem}

Let $\Delta$ be the standard comultiplication in the universal enveloping algebra.

\begin{lem}
\label{comult}

$\Delta PfF_I=\sum_{I'\sqcup I''=I}(-1)^{(I'I'')}PfF_{I'}\otimes
PfF_{I''}$

Here $(-1)^{(I'I'')}$ is a sign of a permutation of the set
$I=\{i_1,...,i_k\}$ that places first the subset $I'\subset I$ and
then places the subset $I''\subset I$.
\end{lem}

\subsection{The action of Pfaffians on tensors}
Let us decribe the action of pfaffians in th tensor representations.

\begin{prop}\label{sta5}

On the base  vectors $e_{-2},e_{-1},e_0,e_1,e_2$ of  the standard
representation of $\mathfrak{o}_{5}$ the pfaffians $PfF_I$ where
$|I|=4$ act as zero operators.
\end{prop}
\proof

The proposition is proved  by direct calculation using the
formulaes, where $a\star b=\frac{1}{2}(ab+ba)$

$PfF_{\widehat{-2}}=F_{0-1}\star F_{-21}-F_{-1-1}\star
F_{-20}+F_{-2-1} \star F_{-10}$

$PfF_{\widehat{-1}}=F_{0-2}\star F_{-21}-F_{-1-2}\star
F_{-20}+F_{-2-2} \star F_{-10}$

$PfF_{\hat{0}}=F_{1-2}\star F_{-21}-F_{-1-2}\star
F_{-2-1}+F_{-2-2}\star F_{-1-1}$

$PfF_{\hat{1}}=F_{1-2}\star F_{-20}-F_{0-2}\star
F_{-2-1}+F_{-2-2}\star F_{0-1}$

$PfF_{\hat{2}}=F_{1-2}\star F_{-10}-F_{0-2}\star
F_{-1-1}+F_{-1-2}\star F_{0-1}$

\endproof

Prove an analog of the previous statement in an arbitrary dimension

\begin{prop}
\label{stan} On the base  vectors $e_{-n},...,e_n$ of the standard
representation of $\mathfrak{o}_{N}$ the pfaffians $PfF_I$ for
$|I|>2$ act as zero operators.
\end{prop}

Put $q=4$, $p=k-4$ in Lemma \ref{minors1}. One has

$$PfF_Ie_j=\sum_{I'\sqcup I''=I,
|I'|=k-4,|I''|=4}\frac{(\frac{p}{2})!(\frac{q}{2})!}{(\frac{k}{2})!}(-1)^{(I'I'')}PfF_{I'}PfF_{I''}e_j.$$

If $j\notin I''$, then obviously  $PfF_{I''}e_j=0$. If  $j\in I''$,
then  using Proposition \ref{sta5} one also obtains
$PfF_{I''}e_j=0$.

\endproof

Let us find an action of a pfaffian of the order  $k$ on a tensor
product of  $<\frac{k}{2}$ vectors,  that is on a tensor product
$e_{r_2}\otimes e_{r_4}\otimes ...\otimes e_{r_t}$, where $t<k$.

\begin{prop}
\label{er} $PfF_Ie_{r_2}\otimes e_{r_4}...\otimes e_{r_t}=0$ where
$t<k$
\end{prop}

\proof The following formulae takes place $\Delta
PfF_I=\sum_{I'\sqcup I''=I}(-1)^{(I'I'')}PfF_{I'}\otimes PfF_{I''}$
(Lemma \ref{comult}).

By definition one has $PfF_Ie_{r_2}\otimes e_{r_4}\otimes...\otimes
e_{r_k}=(\Delta^kPfF_I)e_{r_2}\otimes e_{r_4}\otimes...\otimes
e_{r_k}$. Since $t<k$, the comultiplication $\Delta^kPfF_I$ contains
only  summands in which on some place the pfaffian stands whose
indexing set $I$ satisfies  $|I|\geq 4$ (Lemma \ref{comult}). From
Proposition \ref{sta5} it follows that every such a summand acts as
a zero operator. \endproof

Find an action of a pfaffian of the order  $k$ on a tensor product
of $\frac{k}{2}$ vector,  that is on the tensor product
$e_{r_2}\otimes e_{r_4}\otimes ...\otimes e_{r_k}$.

\begin{prop}
\label{Pfk}If $\{r_2,r_4...,r_k\}$ is not contained in $ I$, then
$PfF_Ie_{r_2}\otimes e_{r_4}\otimes...\otimes e_{r_k}=0$.

 Otherwise
take a permutation $\gamma$ of $I$, such that
$(\gamma(i_1),\gamma(i_2),...,\gamma(i_k))=(r_1,r_2,r_3,...,r_{k-1},r_k)$.
Then \begin{center}$PfF_Ie_{r_2}\otimes...\otimes
e_{r_k}=(-1)^{\gamma}(-1)^{\frac{k(k-1)}{2}}\sum_{\delta\in
Aut(r_1,r_3,...,r_{k-1})}(-1)^{\delta}e_{-\delta(r_1)}\otimes
e_{-\delta(r_3)}\otimes ...\otimes e_{-\delta(r_{k-3})}\otimes
e_{-\delta(r_{k-1})}.$\end{center}

\end{prop}

\proof By definition one has $$PfF_Ie_{r_2}\otimes
e_{r_4}\otimes...\otimes e_{r_k}=(\Delta^kPfF_I)e_{r_2}\otimes
e_{r_4}\otimes...\otimes e_{r_k}.$$ Applying many times the formulae
for comultiplication one obtains
$$\Delta^{\frac{k}{2}}PfF_I=\sum_{I^1\sqcup...\sqcup
I^k}(-1)^{(I^1...I^k)}PfF_{I^1}\otimes...\otimes PfF_{I^k}.$$ Using
Proposition \ref{er} one gets  that, only the summands for which
$|I^j|=2, j=1,...,k$ are nonzero operators.



Hence the summation over divisions  can be written in the following
manner.

\begin{center}$PfF_Ie_{r_2}\otimes e_{r_4}\otimes...\otimes
e_{r_k}=\frac{1}{2^{\frac{k}{2}}}\sum_{\sigma\in
S_k}(-1)^{\sigma}F_{-\sigma(i_1)\sigma(i_2)}\otimes...\otimes
F_{-\sigma(i_{k-1})\sigma(i_k)}(e_{r_2}\otimes...\otimes
e_{r_k})=\frac{1}{2^{\frac{k}{2}}}\sum_{\sigma\in
S_k}(-1)^{\sigma}F_{-\sigma(i_1)\sigma(i_2)}e_{r_2}\otimes...\otimes
F_{-\sigma(i_{k-1})\sigma(i_k)}e_{r_k}.$\end{center}

Consider the expression $F_{-\sigma(i_1)\sigma(i_2)}e_{r_2}$. This
is $e_{-\sigma(i_1)}$ if $\sigma(i_2)=r_2$, this is
$-e_{-\sigma(i_2)}$ if $\sigma(i_1)=r_2$ and zero otherwise. Thus
the summand is nonzero only if the permutation $\sigma$ satisfies
the following condition. In each pair
$(\sigma(i_{2t-1}),\sigma(i_{2t}))$ either $\sigma(i_{2t-1})=r_{2t}$
or   $\sigma(i_{2t})=r_{2t}$.

Show that one can consider only the permutations $\sigma$ such that
$\sigma(i_{2t})=r_{2t}$, that is the permutations of type
$(\sigma(i_1),r_2,\sigma(i_2),r_3...,\sigma(i_{k-1}),r_k)$. But when
only  summands corresponding to such permutations are considered one
must multiply the resulting sum on $2^{\frac{k}{2}}$.

It is enough to prove that the permutations
$\sigma=(\sigma(i_1),\sigma(i_2)=r_2,\sigma(i_3)...,\sigma(r_k))$
and
$\sigma'=(\sigma(i_2)=r_2,\sigma(i_1),\sigma(i_3)...,\sigma(r_k))$
give the same input.

Remind that the input for $\sigma$ is
$$(-1)^{\sigma}F_{-\sigma(i_1)\sigma(i_2)}e_{r_2}\otimes...\otimes
F_{-\sigma(i_{k-1})\sigma(i_k)}e_{r_k}.$$ One has from one hand that
$F_{-\sigma(i_1)\sigma(i_2)}e_{r_2}=e_{-\sigma(i_1)}$ and from the
other hand
$F_{-\sigma'(i_1)\sigma'(i_2)}e_{r_2}=-e_{-\sigma'(i_2)}=-e_{-\sigma(i_1)}$,
 Also one has $(-1)^{\sigma}=-(-1)^{\sigma'}$.
Thus the inputs corresponding to $\sigma$ and $\sigma'$  are the
same.

Hence one can consider the only the permutations  $\sigma$ of type
$(\sigma(i_1),r_2,\sigma(i_2),r_3...,\sigma(i_{k-1}),r_k)$ but
multiplying the resulting sum on $2^{\frac{k}{2}}$.

For the permutation  $\sigma$ of type
$(\sigma(i_1),r_2,\sigma(i_2),r_3...,\sigma(i_{k-1}),r_k)$  using
the definition of $\gamma$ one gets

\begin{center}$(-1)^{\sigma}F_{-\sigma(i_1)\sigma(i_2)}e_{r_2}\otimes...\otimes
F_{-\sigma(i_{k-1})\sigma(i_k)}e_{r_k}=
(-1)^{(\sigma(i_1)r_2,...,\sigma(i_{k-1})r_k)}e_{-\sigma(i_1)}\otimes
e_{-\sigma(i_3)}\otimes...\otimes
e_{-\sigma(i_k)}=(-1)^{\frac{k(k-1)}{2}}(-1)^{\gamma}(-1)^{\delta}e_{-\delta(r_1)}\otimes
...\otimes e_{-\delta(r_{k-1})}.$\end{center} Here $\delta$ is a
permutation of the set $\{r_1,r_3,...,r_{k-3},r_k\}$.

The equality
$(-1)^{\frac{k(k-1)}{2}}(-1)^{\delta}(-1)^{\gamma}=(-1)^{\sigma}$
was used.


 Taking the summation over all permutations $\delta$,  one gets
 $$PfF_Ie_{r_2}\otimes e_{r_4}\otimes...\otimes e_{r_k}=(-1)^{\frac{k(k-1)}{2}}(-1)^{\gamma}\sum_{\delta\in
Aut(r_1,...,r_{k-1})}(-1)^{\delta}e_{-\delta(r_1)}\otimes...\otimes
e_{-\delta(r_{k-1})}. $$ \endproof

Finally from the formula $PfF_Ie_{r_2}\otimes e_{r_4}\otimes
...\otimes e_{r_t}=(\Delta^tPfF_I)e_{r_2}\otimes
e_{r_4}\otimes...\otimes e_{r_t}$, as in the proof of Proposition
\ref{Pfk}, one gets the formulae of the action on an arbitrary
tensor $e_{r_2}\otimes...\otimes e_{r_t}$.

\begin{prop}
\label{Pft} $PfF_Ie_{r_2}\otimes e_{r_4}...\otimes
e_{r_t}=\sum_{\{j_2,j_4...,j_k\}\subset\{2,4...,t\}}Pf^{j_2,j_4...,j_t}F_Ie_{r_2}\otimes
e_{r_4}\otimes...\otimes e_{r_t}$. Here $Pf^{j_2,j_4...,j_k}F_I$
acts on the tensor multiples with numbers $j_2,j_4...,j_k$. It's
action is described by  Proposition \ref{Pfk}

\end{prop}

\subsection{Pfaffians and representation of the algebra $\mathcal{T}$}

Let us give a relation between the representation of the Lie algebra
$\mathcal{T}$ defined by the formula \ref{repr} through the action
of noncommutative pfaffians.

As a corollary of Proposition \ref{Pft} we get the following
proposition

\begin{prop}
\label{corte}In the case $\mathfrak{o}_{2n+1}$ if $I\subset J$,
$|J|=2n$ and $v=e_{i_1}\otimes...\otimes e_{i_{n}}$, then
$PfF_{I}PfF_Jv=C PfF_Jv$, where the constant $C$ depend only on $n$.
\end{prop}

Now let $N=2n+1$ and $|I|=2n$. Then a subset $I$ which consists of
$2n$ elements is of type $\{1,...,N\}\setminus i$. Put
$$PfF_I:=PfF_{\widehat{i}}.$$ Then $PfF_{PfF_{I''}J}=0$ for $I''>2$.

$$[PfF_{\widehat{ i}},PfF_{\widehat{j}}]=\frac{1}{n}\sum_{k\neq i,j} (-1)^{(I'I'')}PfF_{\widehat{ikj}}PfF_{\widehat{k}}$$

The sign $(-1)^{(I'I'')}$ is defined as follows. For an index $s$
denote as $\overline{s}$ either $s$ for $s<i$, оr $s-1$ for  $s>i$.

Then $(-1)^{(I'I'')}$  equals
$(-1)^{n-\overline{j}+(n-1)-\overline{k}}=(-1)^{\overline{j}+\overline{k}-1}$
  in the case $j<k$  and $(-1)^{\overline{j}+\overline{k}}$   in the
  case
  $j>k$. Denote this sign as $s_{jk}$. One has
  $s_{jk}=-s_{kj}$.

Using the theorem \ref{corte} one obtains that for a vector  $v$
from a representation of $\mathfrak{o}_{2n+1}$ with the highest
weight  $(1,...,1)$ the following holds.

\begin{lem}\label{prop8}
$$[PfF_{\widehat{ i}},PfF_{\widehat{j}}]v=C\sum_{k\neq i,j}s_{jk}
PfF_{\widehat{k}}v,$$ where $C$ is some constant.

\end{lem}

The following theorem is proved.

\begin{thm}
\label{rept} The representation of the algebra $\mathcal{T}$ on the
space of skew-symmetric tensors with $n$ indices, given by the
formula \eqref{repr} is given also by the formula

$$e_{a_1}\wedge...\wedge e_{b_n}\mapsto \frac{1}{\sqrt{C}}PfF_{\{a_1,...,b_n\}},$$ where the constant is taken from Lemma
\ref{prop8} and the Pfaffian is considered as an operator acting on
the space of skew-symmetric tensors with $n$ indices.
\end{thm}

\section{Knizhnik-Zamolodchikov equations and Capelli elements}
\label{K}

Now let us give an interpretation of the Kniznik-Zamolodchikov type
equations associated with the root system $B$ constructed for the
algebra $\mathcal{T}$ and it's representation \ref{repr} as a
Kniznik-Zamolodchikov type equations constructed for  higher order
Capelli central elements.

This fact is an immediate corollary  of Theorem \ref{rept}.

Introduce elements.

\begin{equation}\label{cladip}\tau_{ij}=\sum_{I,|I|=2n} 1\otimes...\otimes \rho(PfF_I)\otimes ...\otimes \rho(PfF_{-I})\otimes...1,\end{equation}

where the pfaffians occur on places $i,j$, and $\rho$ is a
representation of $\mathfrak{o}_{2n+1}$ on the space of
skew-symmetric tensors with $n$ indices,

 \begin{equation}\label{clamp}\mu_{ij}=\sum_s 1\otimes...\otimes \rho(Pf F_I)\otimes ...\otimes \rho(Pf F_I)\otimes...1,\end{equation}

\begin{equation}\label{clanp}\nu_{i}=\frac{1}{2}\sum_s 1\otimes.. ...\otimes \rho(PfF_IPfF_{-I}+PfF_IPfF_I)\otimes...1.\end{equation}

\begin{thm}
The action of  elements \ref{cladi} -\ref{clan} and
\ref{cladip}-\ref{clanp} on the space of skew-symmetric tensors
with $n$ indices coincide.
\end{thm}

As a corollary we get

\begin{thm}
The elements \ref{cladip} -\ref{clanp} satisfy the commutation
relation for the coefficients of the Kniznik-Zamolodchikov type
equations associated with the root system $B$
\end{thm}

\end{document}